\documentclass[10pt]{article}
\usepackage{amsmath,amssymb,amsfonts,latexsym,verbatim,bbm,mathrsfs}
\newcommand{\Proof}{$\underline{\it Proof:}$\\}
\newtheorem{lemma}{Lemma}
\newtheorem{theo}{Theorem}

\newtheorem{coro}{Corollary}
\def\square{\mathbin{\vrule height1ex width.5em}}

\begin{document}

\begin{center}
  \huge{\bf Asymptotic results on the moments of the ratio of the random sum of squares
  to the square of the random sum}
\end{center}

\vspace{0.5cm}

\begin{center}
\large{\bf Sophie A. Ladoucette \footnote{Katholieke Universiteit
Leuven, Department of Mathematics, W. de Croylaan 54, B-3001
Leuven, Belgium. Email: sophie.ladoucette@wis.kuleuven.be}}
\end{center}

\vspace{0.7cm}

\begin{center}
  {\bf April 7, 2005}
\end{center}

\vspace{1cm}

\textbf{Abstract:} Let $\{X_1, X_2, \ldots\}$ be a sequence of
positive independent and identically distributed random variables
of Pareto-type with index $\alpha>0$ and let $\{N(t); \, t\geq
0\}$ be a mixed Poisson process independent of the $X_i$'s. For
$t\geq 0$, define:
\begin{equation*}
  T_{N(t)}:=\frac{X_1^2 + X_2^2 + \cdots + X_{N(t)}^2}
  {\left(X_1 + X_2 + \cdots + X_{N(t)}\right)^2}
\end{equation*}
if $N(t)\geq 1$ and $T_{N(t)}:=0$ otherwise.\\

We derive the limiting behavior of the $k$th moment of $T_{N(t)}$,
$k\in\mathbb{N}$, by using the theory of functions of regular
variation and an integral representation for
$\mathbb{E}\left\{T_{N(t)}^k\right\}$. We also point out the
connection between $T_{N(t)}$ and the sample coefficient of
variation which is a popular risk measure in practical
applications.

\vspace{0.5cm}


\vspace{0.5cm}


\section{Introduction}

Let $\{X_1, X_2, \ldots\}$ be a sequence of independent and
identically distributed (i.i.d.) positive random variables with
distribution function $F$ and let $\{N(t); \, t\geq 0\}$ be a
counting process independent of the $X_i$'s. For $t\geq 0$,
define:
\begin{equation*}
  T_{N(t)}:=\frac{X_1^2 + X_2^2 + \cdots + X_{N(t)}^2}
  {\left(X_1 + X_2 + \cdots + X_{N(t)}\right)^2}
\end{equation*}
if $N(t)\geq 1$ and $T_{N(t)}:=0$ otherwise.\\

Denote by $T_n$ the random variable $T_{N(t)}$ when the counting
process $\{N(t); \, t\geq 0\}$ is non-random. An asymptotic
analysis of $T_n$ is provided by Albrecher and Teugels
\cite{at04}, assuming the distribution function $F$ of $X_1$ to be
of \textit{Pareto-type} with positive index $\alpha$. In
particular, they determine limits in distribution for the properly
normalized quantity $T_n$ and derive the limiting behavior of
arbitrary moments of $T_n$, generalizing earlier results
pertaining to $\mathbb{E}\,T_n$ by Fuchs et al.
\cite{fjt01}.\\

In this paper, we focus on \textit{moment convergence} in deriving
the asymptotic behavior of the $k$th moment of $T_{N(t)}$,
$k\in\mathbb{N}$, under the \textit{mixed Poisson} assumption for
the counting process $\{N(t); \, t\geq 0\}$. The distribution
function $F$ of $X_1$ is still assumed to be of Pareto-type with
positive index $\alpha$. The appropriate definitions are recalled
in Section \ref{prel} along with some properties which prove to be
useful later on. The results of the paper are obtained by using
the theory of functions of \textit{regular variation} (e.g.
Bingham et al. \cite{BGT}) and an integral representation for
$\mathbb{E}\left\{T_{N(t)}^k\right\}$ in terms of the probability
generating function of $N(t)$ and the \textit{Laplace transform}
of $X_1$, following in that the basis for the analysis in
Albrecher and Teugels \cite{at04}.\\

When $X_1$ has a Pareto-type distribution function $F$ with
positive index $\alpha$, its moment of order $\beta>0$ is:
\begin{equation*}
  \mu_\beta:=\mathbb{E}\left\{X_1^\beta\right\}=\beta\int_0^\infty
  x^{\beta-1}\left(1-F(x)\right) \, dx \leq\infty
\end{equation*}
which is finite if $\beta<\alpha$ but infinite whenever
$\beta>\alpha$.\\

As pointed out by Albrecher and Teugels \cite{at04}, both the
numerator and the denominator defining $T_{N(t)}$ exhibit an
erratic behavior if $\mu_1=\infty$, whereas this is the case only
for the numerator if $\mu_1<\infty$ and $\mu_2=\infty$. When $X_1$
has a Pareto-type distribution function $F$ with positive index
$\alpha$, then $\mu_1<\infty$ if $\alpha>1$ while $\mu_2<\infty$
as soon as $\alpha>2$. Since the asymptotic behavior of $T_{N(t)}$
is influenced by the finiteness of $\mu_1$ and/or $\mu_2$,
different kinds of results show up depending on the range of
$\alpha$ as presented in Section \ref{sec_results}. In Section
\ref{ccl}, we give some concluding remarks. Incidentally, we point
out the link existing between $T_{N(t)}$ and the \textit{sample
coefficient of variation} of a random sample $X_1, \ldots,
X_{N(t)}$ from a positive random variable $X$ of sample size
$N(t)$ from an integer valued distribution. In a forthcoming
paper, we will take advantage from this link to derive asymptotic
properties of the sample coefficient of variation.

\section{Preliminaries}\label{prel}

Recall that a counting process $\{N(t); \, t\geq 0\}$ is called a
mixed Poisson process if $\{N(t); \, t\geq 0\}=\{\tilde{N}(\Lambda
t); \ t\geq 0\}$, where the mixing random variable $\Lambda$ is
positive and $\{\tilde{N}(t); \ t\geq 0\}$ is a homogeneous
Poisson process with intensity $1$ independent of $\Lambda$. For
each fixed $t\geq 0$, the random variable $N(t)$ has a mixed
Poisson distribution, with mixing distribution the distribution
function $H$ of $\Lambda$, given by:
\begin{equation*}
p_n(t):=\mathbb{P}[N(t)=n]=\mathbb{E}\left\{\frac{(\Lambda
t)^n}{n!} \, e^{-\Lambda t}\right\}=\int_0^\infty \frac{(\lambda
t)^n}{n!} \, e^{-\lambda t} \, dH(\lambda), \quad n\in\mathbb{N}.
\end{equation*}

If the distribution function $H$ is degenerate at a single point
$\lambda \in (0,\infty)$ then the counting process is a
homogeneous Poisson process with intensity $\lambda$. The latter
plays a crucial role in practical applications. In particular, it
is the most popular among all claim number processes in the
actuarial literature. Also, the mixed Poisson process, introduced
to actuaries by Dubourdieu \cite{du38}, has always been very
popular among (re)insurance modelers. It has found many
applications in (re)insurance mathematics because of its
flexibility, its success in actuarial data fitting and its
property of being more dispersed than the Poisson process. For a
general overview on mixed Poisson processes, we refer to the
monograph by Grandell \cite{gr97}.\\

Another way of highlighting the role of the random variable
$\Lambda$ can be expressed by the observation that:
\begin{equation*}
  \frac{N(t)}{t}\stackrel{a.s.}{\longrightarrow}\Lambda \quad \text{as} \quad t \rightarrow \infty
\end{equation*}
where $\stackrel{a.s.}{\longrightarrow}$ stands for almost sure
convergence.\\

For a fixed time $t\geq 0$, the probability generating function of
$N(t)$ is denoted by $Q_t(.)$ and satisfies:
\begin{equation*}
  Q_t(z):=\mathbb{E}\left\{z^{N(t)}\right\}=\sum_{n=0}^\infty p_n(t) \,
  z^n=\mathbb{E}\left\{e^{-t(1-z)\Lambda}\right\}, \quad |z|\leq 1.
\end{equation*}

The $r$th derivative of $Q_t(z)$ with respect to $z$ is denoted by
$Q_t^{(r)}(z)$ and is defined for $|z|<1$. It can be expressed in
terms of expectations as:
\begin{equation*}
  Q_t^{(r)}(z)=r! \, \mathbb{E}\left\{\binom{N(t)}{r} \,
  z^{N(t)-r}\right\}=t^r \, \mathbb{E}\left\{e^{-t(1-z)\Lambda}\Lambda^r\right\}.
\end{equation*}

We define the auxiliary quantities
$q_r(w):=\mathbb{E}\left\{e^{-w\Lambda}\Lambda^r\right\}$, $r
\in\mathbb{N}$, $w\geq 0$, where
$q_r(0)=\mathbb{E}\left\{\Lambda^r\right\}\leq\infty$. Notice that
for all $0\leq w<t$ and $r\in\mathbb{N}$, the following identity
holds:
\begin{equation*}
\frac{1}{t^r}Q_t^{(r)}\left(1-\frac{w}{t}\right)=q_r(w)
\end{equation*}
where the right-hand side does no longer depend on $t$.\\

Before giving an easy but useful result on the moment condition
for $\Lambda$, note that for any $\beta>0$:
\begin{equation}\label{intbeta}
  \int_0^\infty w^{\beta-1} \, q_r(w) \, dw
  = \Gamma(\beta) \, \mathbb{E}\left\{\Lambda^{r-\beta}\right\}
\end{equation}
where $\Gamma(.)$ denotes the gamma function.

\vspace{0.5cm}

\begin{lemma}\label{lemmamomcond}
  Let $\Lambda$ be a positive random variable with
  distribution function $H$. Then for all $0<r\leq s$,
  $\mathbb{E}\left\{\Lambda^s\right\}<\infty\Rightarrow\mathbb{E}\left\{\Lambda^r\right\}<\infty$
  and
  $\mathbb{E}\left\{\Lambda^{-s}\right\}<\infty\Rightarrow\mathbb{E}\left\{\Lambda^{-r}\right\}<\infty$.
\end{lemma}

\vspace{0.5cm}

\Proof Let $0<r\leq s$. Assume that
$\mathbb{E}\left\{\Lambda^s\right\}<\infty$. Then:
\begin{equation*}
  \mathbb{E}\left\{\Lambda^r\right\}=\int_0^1 \lambda^r \, dH(\lambda)
  +\int_1^\infty \lambda^r \, dH(\lambda) \leq \mathbb{P}[\Lambda\leq
  1]+\mathbb{E}\left\{\Lambda^s\right\}<\infty.
\end{equation*}

Now, assume that $\mathbb{E}\left\{\Lambda^{-s}\right\}<\infty$.
Then:
\begin{equation*}
  \mathbb{E}\left\{\Lambda^{-r}\right\}=\int_0^1 \lambda^{-r} \, dH(\lambda)
  +\int_1^\infty \lambda^{-r} \, dH(\lambda) \leq \mathbb{E}\left\{\Lambda^{-s}\right\}
  +\mathbb{P}[\Lambda>1]<\infty. \ \square
\end{equation*}

\vspace{0.5cm}

Recall that the process $\{X_i; \, i\geq 1\}$ consists of i.i.d.
positive random variables with distribution function $F$. As
specified above, our asymptotic results are derived under the
condition that $F$ is of Pareto-type with positive index $\alpha$,
or equivalently that $F$ has a regularly varying tail at $\infty$
with negative index $-\alpha$. This means that the tail of $F$
satisfies:
\begin{equation}\label{defPT}
  1-F(x)\sim x^{-\alpha} \, \ell(x) \quad \text{as} \quad x\rightarrow\infty
\end{equation}
where $\alpha>0$ and $\ell$ is slowly varying at $\infty$.\\

Recall that a measurable and ultimately positive function $f$ on
$\mathbb{R}_+$ is regularly varying at $\infty$ with index
$\gamma\in\mathbb{R}$ (written $f\in \mathrm{RV}_\gamma^\infty$)
if for all $x>0$, $\lim_{t\rightarrow\infty}f(tx)/f(t)=x^\gamma$.
When $\gamma=0$, $f$ is said to be slowly varying at $\infty$.
Similarly, a function $g$ on $\mathbb{R}_+$ is regularly varying
at $0$ with index $\gamma\in\mathbb{R}$ (written $g\in
\mathrm{RV}_\gamma^0$) if for all $x>0$, $\lim_{s\rightarrow
0}g(sx)/g(s)=x^\gamma$. When $\gamma=0$, $g$ is said to be slowly
varying at $0$. For a textbook treatment on the theory of
functions of regular variation, we refer to Bingham et al.
\cite{BGT}.\\

It is well-known that condition (\ref{defPT}) appears as the
essential condition in the extremal domain of attraction problem
of extreme value theory. For a recent treatment, see Beirlant et
al. \cite{bgst04}. When $\alpha\in(0,2)$, the condition
(\ref{defPT}) is also necessary and sufficient for $F$ to belong
to the additive domain of attraction of a non-normal stable law
with exponent $\alpha$ (e.g. Theorem 8.3.1 of Bingham et al.
\cite{BGT}).\\

The common Laplace transform of the $X_i$'s is defined and denoted
by:
\begin{equation*}
  \varphi(s):=\mathbb{E}\left\{e^{-sX_1}\right\}=\int_0^\infty e^{-sx} \, dF(x), \quad s\geq 0.
\end{equation*}

We denote by $\varphi^{(n)}(s)$ the $n$th derivative of
$\varphi(s)$ with respect to $s$. By Lemma 3.1 of Albrecher and
Teugels \cite{at04} and Bingham-Doney's lemma (e.g. Theorem 8.1.6
of Bingham et al. \cite{BGT}), the asymptotic behavior of
$\varphi^{(n)}$ at the origin when $F$ satisfies (\ref{defPT}) is
the following.

\vspace{0.5cm}

\begin{lemma}\label{lemmaderiv}
  If the distribution function $F$ of $X_1$ satisfies
  $1-F(x)\sim x^{-\alpha}\ell(x)$ as $x\rightarrow\infty$ for
  some $\ell\in\mathrm{RV}_0^\infty$ and $\alpha>0$, then:
  \begin{equation*}
  (-1)^n \, \varphi^{(n)}(s) \underset{s\downarrow 0}{\sim}
    \begin{cases}
      \alpha \, \Gamma(n-\alpha) \, s^{\alpha-n} \, \ell\left(\frac{1}{s}\right) & \text{if $n>\alpha$}\\
      \alpha \, \tilde{\ell}\left(\frac{1}{s}\right) &
      \text{if $n=\alpha$ and $\mu_n=\infty$}\\
      \mu_n & \text{if $n<\alpha$ or if $n=\alpha$ and $\mu_n<\infty$}
    \end{cases}
  \end{equation*}
  where $\tilde{\ell}(x)=\int_0^x \frac{\ell(u)}{u} \, du \in \mathrm{RV}_0^\infty$.
\end{lemma}

\vspace{0.5cm}

Now, we give our results.

\section{Results}\label{sec_results}

We start by deriving an integral representation for the $k$th
moment of $T_{N(t)}$ under the mixed Poisson assumption for the
counting process $\{N(t); \, t\geq 0\}$. Note that we do not make
any assumption on the distribution function $F$ of $X_1$.

\vspace{0.5cm}

\begin{lemma}
  Let $t\geq 0$ and $k\in\mathbb{N}\backslash\{0\}$ be fixed.
  Assume that $\{N(t); \, t\geq 0\}$ is a mixed Poisson
  process. The $k$th moment of $T_{N(t)}$ is then given by:
  \begin{equation}\label{momT}
     \mathbb{E}\left\{T_{N(t)}^k\right\}=\sum_{r=1}^{k}\sum_{\begin{subarray}{c}
     k_1, \ldots, k_r\geq 1 \\ k_1+\cdots+k_r=k \end{subarray}}
     \frac{k!}{\prod_{i=1}^{r}k_i!}\frac{B_t(k_1, \ldots, k_r)}{(2k-1)! \, r!}
  \end{equation}
  with:
  \begin{equation}\label{Bt}
    B_t(k_1, \ldots, k_r):=\int_0^\infty s^{2k-1}\prod_{i=1}^{r}\varphi^{(2k_i)}(s)
    \, Q_t^{(r)}(\varphi(s)) \, ds.
  \end{equation}
\end{lemma}

\vspace{0.5cm}

\Proof Using Lemma 2.1 of Albrecher and Teugels \cite{at04}, we
easily derive:
\begin{eqnarray*}
  \mathbb{E}\left\{T_{N(t)}^k\right\} &=& \sum_{n=0}^{\infty}
  p_n(t) \, \mathbb{E}\left\{\left. T_{N(t)}^k \right\vert N(t)=n\right\} = \sum_{n=0}^{\infty}
  p_n(t) \, \mathbb{E}\left\{T_n^k\right\} \cr
  &=& \sum_{n=0}^{\infty} p_n(t)\sum_{r=1}^{k}\sum_{\begin{subarray}{c}
  k_1, \ldots, k_r\geq 1 \\ k_1+\cdots+k_r=k \end{subarray}}\frac{k!}{\prod_{i=1}^{r}k_i!}
  \frac{\binom{n}{r}}{(2k-1)!} \int_0^\infty s^{2k-1}
  \prod_{i=1}^{r}\varphi^{(2k_i)}(s) \, \varphi^{n-r}(s) \, ds \cr
  &=& \sum_{r=1}^{k}\sum_{\begin{subarray}{c} k_1, \ldots, k_r\geq 1 \\
  k_1+\cdots+k_r=k \end{subarray}}\frac{k!}{\prod_{i=1}^{r}k_i!} \frac{1}{(2k-1)!}
  \int_0^\infty s^{2k-1} \prod_{i=1}^{r}\varphi^{(2k_i)}(s) \sum_{n=0}^{\infty} p_n(t)
  \binom{n}{r}\varphi^{n-r}(s) \, ds \cr
  &=& \sum_{r=1}^{k}\sum_{\begin{subarray}{c} k_1, \ldots, k_r\geq 1 \\
  k_1+\cdots+k_r=k \end{subarray}}\frac{k!}{\prod_{i=1}^{r}k_i!} \frac{1}{(2k-1)! \, r!}
  \int_0^\infty s^{2k-1} \prod_{i=1}^{r}\varphi^{(2k_i)}(s) \, Q_t^{(r)}(\varphi(s))
  \, ds. \ \square
\end{eqnarray*}

\vspace{0.5cm}

We rewrite $B_t(k_1, \ldots, k_r)$ in a more convenient form which
proves useful later on. Defining $\psi(s):=\varphi^{-1}(1-s)$ and
substituting $s=\psi\left(\frac{w}{t}\right)$ in (\ref{Bt}) lead
to: \small
\begin{eqnarray*}
  & & B_t(k_1, \ldots, k_r) \cr
  & & \quad = \int_0^\infty s^{2k-1}\prod_{i=1}^{r}\varphi^{(2k_i)}(s) \, Q_t^{(r)}(\varphi(s)) \, ds \cr
  & & \quad = \int_0^t \psi^{2k-1}\left(\frac{w}{t}\right)
  \prod_{i=1}^{r}\varphi^{(2k_i)}\left(\psi\left(\frac{w}{t}\right)\right)
  \, Q_t^{(r)}\left(1-\frac{w}{t}\right) \, \frac{d}{dw}\psi\left(\frac{w}{t}\right) \, dw \cr
  & & \quad = -t^{r-1}\frac{\psi^{2k-1}\left(\frac{1}{t}\right)}
  {\varphi^{(1)}\left(\psi\left(\frac{1}{t}\right)\right)}
  \prod_{i=1}^{r}\varphi^{(2k_i)}\left(\psi\left(\frac{1}{t}\right)\right)
  \int_0^t \left(\frac{\psi\left(\frac{w}{t}\right)}{\psi\left(\frac{1}{t}\right)}\right)^{2k-1}
  \prod_{i=1}^{r}\frac{\varphi^{(2k_i)}\left(\psi\left(\frac{w}{t}\right)\right)}
  {\varphi^{(2k_i)}\left(\psi\left(\frac{1}{t}\right)\right)}
  \, \frac{\varphi^{(1)}\left(\psi(\frac{1}{t})\right)}
  {\varphi^{(1)}\left(\psi(\frac{w}{t})\right)} \, \frac{Q_t^{(r)}(1-\frac{w}{t})}{t^r} \, dw \cr
  & & \quad = -t^{r-1}\frac{\psi^{2k-1}\left(\frac{1}{t}\right)}
  {\varphi^{(1)}\left(\psi\left(\frac{1}{t}\right)\right)}
  \prod_{i=1}^{r}\varphi^{(2k_i)}\left(\psi\left(\frac{1}{t}\right)\right)
  \int_0^\infty \left(\frac{\psi\left(\frac{w}{t}\right)}{\psi\left(\frac{1}{t}\right)}\right)^{2k-1}
  \prod_{i=1}^{r}\frac{\varphi^{(2k_i)}\left(\psi\left(\frac{w}{t}\right)\right)}
  {\varphi^{(2k_i)}\left(\psi\left(\frac{1}{t}\right)\right)}
  \, \frac{\varphi^{(1)}\left(\psi\left(\frac{1}{t}\right)\right)}
  {\varphi^{(1)}\left(\psi(\frac{w}{t})\right)} \, q_r(w) \, \mathbbm{1}_{[0,t)}(w) \, dw \cr
  & & \quad = f_t(k_1, \ldots, k_r)\int_0^\infty g_t(w;k_1, \ldots, k_r) \ dw
\end{eqnarray*}
\normalsize with:
\begin{equation*}
  f_t(k_1, \ldots, k_r):=-t^{r-1}\frac{\psi^{2k-1}\left(\frac{1}{t}\right)}
  {\varphi^{(1)}\left(\psi\left(\frac{1}{t}\right)\right)}
  \prod_{i=1}^{r}\varphi^{(2k_i)}\left(\psi\left(\frac{1}{t}\right)\right)
\end{equation*}
and:
\begin{equation*}
  g_t(w;k_1, \ldots, k_r):=\left(\frac{\psi\left(\frac{w}{t}\right)}
  {\psi\left(\frac{1}{t}\right)}\right)^{2k-1}
  \prod_{i=1}^{r}\frac{\varphi^{(2k_i)}\left(\psi\left(\frac{w}{t}\right)\right)}
  {\varphi^{(2k_i)}\left(\psi\left(\frac{1}{t}\right)\right)}
  \, \frac{\varphi^{(1)}\left(\psi\left(\frac{1}{t}\right)\right)}
  {\varphi^{(1)}\left(\psi(\frac{w}{t})\right)} \, q_r(w) \, \mathbbm{1}_{[0,t)}(w).
\end{equation*}

\vspace{0.5cm}

From now on, the distribution function $F$ of $X_1$ is assumed to
satisfy (\ref{defPT}) for some $\alpha>0$. Here is the first of
our main results pertaining to moment convergence for $T_{N(t)}$.
It concerns the case $\alpha\in(0,1)$.

\vspace{0.5cm}

\begin{theo}
  Assume that $X_1$ is of Pareto-type with index $\alpha\in(0,1)$. Let $\{N(t); \, t\geq 0\}$
  be a mixed Poisson process with mixing random variable $\Lambda$. If
  $\mathbb{E}\left\{\Lambda^{\epsilon}\right\}<\infty$ and
  $\mathbb{E}\left\{\Lambda^{-\epsilon}\right\}<\infty$ for some $\epsilon>0$, then
  for any fixed $k\in\mathbb{N}\backslash\{0\}$:
  \begin{equation*}
    \lim_{t\rightarrow\infty}\mathbb{E}\left\{T_{N(t)}^k\right\}=\frac{k!}{(2k-1)!}
    \sum_{r=1}^{k}\frac{\alpha^{r-1}}{r \, \Gamma^r(1-\alpha)} \, G(r,k)
  \end{equation*}
  where $G(r,k)$ is the coefficient of $x^k$ in the polynomial
  $\left(\sum_{i=1}^{k-r+1}\frac{\Gamma(2i-\alpha)}{i!} \, x^i\right)^r$.
\end{theo}

\vspace{0.5cm}

\Proof Let $k\in\mathbb{N}\backslash\{0\}$ and $\alpha\in(0,1)$ be
fixed. Since
$1-F(x)\underset{x\uparrow\infty}{\sim}x^{-\alpha}\ell(x)$, it
follows from Corollary 8.1.7 in \cite{BGT} that
$1-\varphi(s)\underset{s\downarrow 0}{\sim}\Gamma(1-\alpha) \,
s^{\alpha}\ell\left(\frac{1}{s}\right)$. Hence, we easily deduce
$s=1-\varphi(\psi(s))\underset{s\downarrow
0}{\sim}\Gamma(1-\alpha) \, \psi^\alpha(s) \,
\ell\left(1/\psi(s)\right)$ which leads to:
\begin{equation}\label{rel1}
  \lim_{s\rightarrow 0} s^{-1} \, \psi^\alpha(s) \,
  \ell\left(\frac{1}{\psi(s)}\right)=\frac{1}{\Gamma(1-\alpha)}.
\end{equation}

Relation (\ref{rel1}) learns us that $\psi$ is regularly varying
at 0 with index $1/\alpha$. For $n>\alpha$, we consequently get
$\varphi^{(n)}\circ\psi\in \mathrm{RV}_{1-\frac{n}{\alpha}}^0$
since $\varphi^{(n)}\in \mathrm{RV}_{\alpha-n}^0$ by Lemma
\ref{lemmaderiv}, $\psi\in \mathrm{RV}_{1/\alpha}^0$ and
$\lim_{s\rightarrow 0}\psi(s)=0$.\\

Using Potter's theorem (e.g. Theorem 1.5.6 of Bingham et al.
\cite{BGT}), we therefore obtain the following upper bound for the
integrand in $B_t(k_1, \ldots, k_r)$. Set
$\delta_r:=\frac{\zeta}{2k+r}$ with $\zeta=\epsilon$ if
$\epsilon\in(0,1)$ or $\zeta\in(0,1)$ otherwise. For this chosen
$\delta_r>0$, there exists $C_r=C_r(\delta_r)>1$ such that for all
$t>0$:
\begin{eqnarray*}
  g_t(w;k_1, \ldots, k_r) &\leq& C_r \, w^{\frac{2k-1}{\alpha}}\left(\max\left\{w^{\delta_r},
  w^{-\delta_r}\right\}\right)^{2k-1} w^{\frac{r\alpha-2k}{\alpha}}\left(\max\left\{w^{\delta_r},
  w^{-\delta_r}\right\}\right)^{r+1}w^{\frac{1-\alpha}{\alpha}} q_r(w) \cr
  &=& C_r \, w^{r-1} \, \max\left\{w^{\zeta}, w^{-\zeta}\right\} \, q_r(w) =:h(w).
\end{eqnarray*}

Now, $\int_0^\infty h(w) \, dw<\infty$ if and only if $\int_0^1
w^{r-1-\zeta} \, q_r(w) \, dw<\infty$ and $\int_1^\infty
w^{r-1+\zeta} \, q_r(w) \, dw<\infty$.\\

Since $\zeta\in(0,1)$ and $\zeta\leq \epsilon$, we use
(\ref{intbeta}) together with Lemma \ref{lemmamomcond} to get:
\begin{equation*}
  \int_0^1 w^{r-1-\zeta} \, q_r(w) \, dw \leq \int_0^\infty w^{r-1-\zeta} \, q_r(w) \, dw
  =\Gamma(r-\zeta) \, \mathbb{E}\left\{\Lambda^{\zeta}\right\} < \infty
\end{equation*}
and:
\begin{equation*}
  \int_1^\infty w^{r-1+\zeta} \, q_r(w) \, dw \leq \int_0^\infty w^{r-1+\zeta} \, q_r(w) \, dw
  =\Gamma(r+\zeta) \, \mathbb{E}\left\{\Lambda^{-\zeta}\right\} < \infty.
\end{equation*}

Hence, the function $h$ is integrable.\\

Finally, $\lim_{t\rightarrow\infty}g_t(w;k_1, \ldots, k_r)=
w^{\frac{2k-1}{\alpha}}w^{r-\frac{2k}{\alpha}}w^{\frac{1}{\alpha}-1}q_r(w)=w^{r-1}q_r(w)$.
Thus, applying Lebesgue's theorem on dominated convergence and
using (\ref{intbeta}), we deduce:
\begin{equation*}
  \lim_{t\rightarrow\infty}\int_0^\infty g_t(w;k_1, \ldots, k_r) \, dw
  =\int_0^\infty w^{r-1} \, q_r(w) \, dw=(r-1)!.
\end{equation*}

Using Lemma \ref{lemmaderiv} and relation (\ref{rel1}), we get:
\begin{eqnarray*}
  f_t(k_1, \ldots, k_r) &\underset{t\uparrow\infty}{\sim}& \frac{-t^{r-1}
  \alpha^r\prod_{i=1}^{r}\Gamma(2k_i-\alpha) \, \psi^{r\alpha-1}\left(\frac{1}{t}\right)
  \ell^r\left(\frac{1}{\psi\left(\frac{1}{t}\right)}\right)}{-\alpha \, \Gamma(1-\alpha)
  \, \psi^{\alpha-1}\left(\frac{1}{t}\right)
  \ell\left(\frac{1}{\psi\left(\frac{1}{t}\right)}\right)} \cr
  &=& \alpha^{r-1}\frac{\prod_{i=1}^{r}\Gamma(2k_i-\alpha)}{\Gamma(1-\alpha)}
  \,  t^{r-1}\psi^{\alpha(r-1)}\left(\frac{1}{t}\right)
  \ell^{r-1}\left(\frac{1}{\psi\left(\frac{1}{t}\right)}\right)
\end{eqnarray*}
so that:
\begin{equation*}
  \lim_{t\rightarrow\infty}f_t(k_1, \ldots, k_r)=
  \frac{\alpha^{r-1}\prod_{i=1}^r\Gamma(2k_i-\alpha)}{\Gamma^r(1-\alpha)}.
\end{equation*}

Therefore, we obtain:
\begin{equation*}
  \lim_{t\rightarrow\infty}B_t(k_1, \ldots, k_r)=
  \frac{(r-1)! \, \alpha^{r-1}\prod_{i=1}^r\Gamma(2k_i-\alpha)}{\Gamma^r(1-\alpha)}.
\end{equation*}

Summing up over all $r\in\{1, \ldots, k\}$ in (\ref{momT}), we
arrive at:
\begin{equation*}
  \lim_{t\rightarrow\infty}\mathbb{E}\left\{T_{N(t)}^k\right\}=
  \frac{k!}{(2k-1)!}\sum_{r=1}^{k}\frac{\alpha^{r-1}}{r \, \Gamma^r(1-\alpha)}
  \sum_{\begin{subarray}{c} k_1, \ldots, k_r\geq 1 \\ k_1+\cdots+k_r=k \end{subarray}}
  \prod_{i=1}^r\frac{\Gamma(2k_i-\alpha)}{k_i!}.
\end{equation*}

Finally, Albrecher and Teugels \cite{at04} have observed that:
\begin{equation*}
  G(r,k):=\sum_{\begin{subarray}{c} k_1, \ldots, k_r\geq 1 \\ k_1+\cdots+k_r=k
  \end{subarray}} \prod_{i=1}^r\frac{\Gamma(2k_i-\alpha)}{k_i!}
\end{equation*}
can be read off as the coefficient of $x^k$ in the $r$-fold
product $\left(\sum_{i=1}^{k-r+1}\frac{\Gamma(2i-\alpha)}{i!} \,
x^i\right)^r$. \ $\square$

\vspace{0.5cm}

Our next result deals with the case $\alpha=1$ and $\mu_1=\infty$.

\vspace{0.5cm}

\begin{theo}\label{theo1}
  Assume that $X_1$ is of Pareto-type with index $\alpha=1$ and that $\mu_1=\infty$. Let
  $\{N(t); \, t\geq 0\}$ be a mixed Poisson process with mixing random variable $\Lambda$.
  If $\mathbb{E}\left\{\Lambda^{\epsilon}\right\}<\infty$
  and $\mathbb{E}\left\{\Lambda^{-\epsilon}\right\}<\infty$
  for some $\epsilon>0$, then for any fixed $k\in\mathbb{N}\backslash\{0\}$:
  \begin{equation*}
    \mathbb{E}\left\{T_{N(t)}^k\right\}\sim \frac{1}{2k-1}
    \frac{\ell(a_t)}{\tilde{\ell}(a_t)}
    \quad \text{as} \quad t\rightarrow\infty
  \end{equation*}
  where $\tilde{\ell}(x)=\int_0^x \frac{\ell(u)}{u} \, du \in\mathrm{RV}_0^\infty$
  and $(a_t)_{t>0}$ is a sequence defined by
  $\lim_{t\rightarrow\infty} t \, a_t^{-1} \, \tilde{\ell}(a_t)=1$.
\end{theo}

\vspace{0.5cm}

\Proof Let $k\in\mathbb{N}\backslash\{0\}$ be fixed. Since
$\mu_1=\infty$ and $1-F(x)\underset{x\uparrow
\infty}{\sim}x^{-1}\ell(x)$ for some $\ell\in
\mathrm{RV}_0^\infty$, it follows by Lemma \ref{lemmaderiv} that
$\varphi^{(1)}(s)\underset{s\downarrow
0}{\sim}-\tilde{\ell}\left(\frac{1}{s}\right)$ and then that
$1-\varphi(s)\underset{s\downarrow 0}{\sim}s \,
\tilde{\ell}\left(\frac{1}{s}\right)$ with
$\tilde{\ell}(x)=\int_0^x \frac{\ell(u)}{u} \, du\in
\mathrm{RV}_0^\infty$. Since
$s=1-\varphi(\psi(s))\underset{s\downarrow
0}{\sim}\psi(s)\tilde{\ell}\left(1/\psi(s)\right)$, we obtain:
\begin{equation}\label{rel2}
  \lim_{s\rightarrow 0} s^{-1} \, \psi(s) \, \tilde{\ell}\left(\frac{1}{\psi(s)}\right)=1.
\end{equation}

Relation (\ref{rel2}) learns us that $\psi\in \mathrm{RV}_1^0$,
leading to $\lim_{s\rightarrow 0}\psi(s)=0$. For $n\geq 2$, we
consequently get $\varphi^{(n)}\circ\psi\in\mathrm{RV}_{1-n}^0$
since $\varphi^{(n)}\in \mathrm{RV}_{1-n}^0$ by Lemma
\ref{lemmaderiv}. Moreover,
$\varphi^{(1)}\circ\psi\in\mathrm{RV}_0^0$ since
$\varphi^{(1)}\in\mathrm{RV}_0^0$.\\

Using Potter's theorem, we therefore obtain the following upper
bound for the integrand in $B_t(k_1, \ldots, k_r)$. Set
$\delta_r:=\frac{\zeta}{2k+r}$ with $\zeta=\epsilon$ if
$\epsilon\in(0,1)$ or $\zeta\in(0,1)$ otherwise. For this chosen
$\delta_r>0$, there exists $C_r=C_r(\delta_r)>1$ such that for all
$t>0$:
\begin{eqnarray*}
  g_t(w;k_1, \ldots, k_r) &\leq& C_r \, w^{2k-1}\left(\max\left\{w^{\delta_r},
  w^{-\delta_r}\right\}\right)^{2k-1} w^{r-2k}\left(\max\left\{w^{\delta_r},
  w^{-\delta_r}\right\}\right)^{r+1} q_r(w) \cr
  &=& C_r \, w^{r-1} \, \max\left\{w^{\zeta}, w^{-\zeta}\right\} \, q_r(w) =:h(w).
\end{eqnarray*}

Now, $\int_0^\infty h(w) \, dw<\infty$ if and only if $\int_0^1
w^{r-1-\zeta} \, q_r(w) \, dw<\infty$ and $\int_1^\infty
w^{r-1+\zeta} \, q_r(w) \, dw<\infty$.\\

Since $\zeta\in(0,1)$ and $\zeta\leq \epsilon$, we use
(\ref{intbeta}) together with Lemma \ref{lemmamomcond} to get:
\begin{equation*}
  \int_0^1 w^{r-1-\zeta} \, q_r(w) \, dw \leq \int_0^\infty w^{r-1-\zeta} \, q_r(w) \, dw
  =\Gamma(r-\zeta) \, \mathbb{E}\left\{\Lambda^{\zeta}\right\} < \infty
\end{equation*}
and:
\begin{equation*}
  \int_1^\infty w^{r-1+\zeta} \, q_r(w) \, dw \leq \int_0^\infty w^{r-1+\zeta} \, q_r(w) \, dw
  =\Gamma(r+\zeta) \, \mathbb{E}\left\{\Lambda^{-\zeta}\right\} < \infty.
\end{equation*}

Hence, the function $h$ is integrable.\\

Finally, $\lim_{t\rightarrow\infty}g_t(w;k_1, \ldots, k_r)=
w^{2k-1}w^{r-2k}q_r(w)=w^{r-1}q_r(w)$. Thus, applying Lebesgue's
theorem on dominated convergence and using (\ref{intbeta}), we
deduce:
\begin{equation*}
  \lim_{t\rightarrow\infty}\int_0^\infty g_t(w;k_1, \ldots, k_r) \, dw
  =\int_0^\infty w^{r-1} \, q_r(w) \, dw=(r-1)!.
\end{equation*}

Using Lemma \ref{lemmaderiv}, relation (\ref{rel2}) and
$\varphi^{(1)}(\psi(s))\underset{s\downarrow
0}{\sim}-\tilde{\ell}\left(1/\psi(s)\right)$, we get:
\begin{eqnarray*}
  f_t(k_1, \ldots, k_r) &\underset{t\uparrow\infty}{\sim}& t^{r-1}
  \prod_{i=1}^{r}\Gamma(2k_i-1) \, \psi^{r-1}\left(\frac{1}{t}\right)
  \frac{\ell^r\left(1/\psi\left(\frac{1}{t}\right)\right)}
  {\tilde{\ell}\left(1/\psi\left(\frac{1}{t}\right)\right)} \cr
  &\underset{t\uparrow\infty}{\sim}& \prod_{i=1}^{r}\Gamma(2k_i-1)
  \left(\frac{\ell\left(1/\psi\left(\frac{1}{t}\right)\right)}
  {\tilde{\ell}\left(1/\psi\left(\frac{1}{t}\right)\right)}\right)^r.
\end{eqnarray*}

Therefore, we obtain:
\begin{equation*}
  B_t(k_1, \ldots, k_r) \underset{t\uparrow\infty}{\sim} (r-1)! \prod_{i=1}^{r}\Gamma(2k_i-1)
  \left(\frac{\ell\left(1/\psi\left(\frac{1}{t}\right)\right)}
  {\tilde{\ell}\left(1/\psi\left(\frac{1}{t}\right)\right)}\right)^r.
\end{equation*}

Since
$\lim_{t\rightarrow\infty}\frac{\ell\left(1/\psi\left(\frac{1}{t}\right)\right)}
{\tilde{\ell}\left(1/\psi\left(\frac{1}{t}\right)\right)}=0$, only
the summand with $r=1$ contributes to the dominating asymptotic
term of (\ref{momT}). Hence, we get:
\begin{equation*}
  \mathbb{E}\left\{T_{N(t)}^k\right\} \sim \frac{1}{2k-1}
  \frac{\ell\left(1/\psi\left(\frac{1}{t}\right)\right)}
  {\tilde{\ell}\left(1/\psi\left(\frac{1}{t}\right)\right)}
  \quad \text{as} \quad t\rightarrow\infty.
\end{equation*}

One easily notes that relation (\ref{rel2}) is equivalent to
$\lim_{t\rightarrow 0} t \, \psi\left(\frac{1}{t}\right) \,
\tilde{\ell}\left(1/\psi\left(\frac{1}{t}\right)\right)=1$.
Defining a sequence $(a_t)_{t>0}$ by $\lim_{t\rightarrow\infty} t
\, a_t^{-1} \, \tilde{\ell}(a_t)=1$ implies that
$\psi\left(\frac{1}{t}\right)\underset{t\uparrow
\infty}{\sim}\frac{1}{a_t}$. By virtue of the uniform convergence
theorem for slowly varying functions (e.g. Theorem 1.2.1 of
Bingham et al. \cite{BGT}), we thus get
$\ell\left(1/\psi\left(\frac{1}{t}\right)\right)\underset{t\uparrow
\infty}{\sim}\ell(a_t)$ and
$\tilde{\ell}\left(1/\psi\left(\frac{1}{t}\right)\right)\underset{t\uparrow
\infty}{\sim}\tilde{\ell}(a_t)$. Consequently, we finally arrive
at:
\begin{equation*}
  \mathbb{E}\left\{T_{N(t)}^k\right\} \sim \frac{1}{2k-1}\frac{\ell(a_t)}{\tilde{\ell}(a_t)}
  \quad \text{as} \quad t\rightarrow\infty. \ \square
\end{equation*}

\vspace{0.5cm}

In the following result, the case $\alpha\in(1,2)$ (including
$\alpha=1$ if $\mu_1<\infty$) is of interest.

\vspace{0.5cm}

\begin{theo}\label{theo12}
  Assume that $X_1$ is of Pareto-type with index $\alpha\in(1,2)$ (including
  $\alpha=1$ if $\mu_1<\infty$). Let $\{N(t); \, t\geq 0\}$ be a mixed Poisson process
  with mixing random variable $\Lambda$. For any fixed $k\in\mathbb{N}\backslash\{0\}$, if
  $\mathbb{E}\left\{\Lambda^{\epsilon}\right\}<\infty$ and
  $\mathbb{E}\left\{\Lambda^{-k(\alpha-1)-\epsilon}\right\}<\infty$ for some $\epsilon>0$, then:
  \begin{equation*}
    \mathbb{E}\left\{T_{N(t)}^k\right\}\sim
    \frac{\alpha}{\mu_1^\alpha} \, B(2k-\alpha,\alpha) \, \mathbb{E}\left\{\Lambda^{1-\alpha}\right\}
    t^{1-\alpha} \, \ell(t) \quad \text{as} \quad t\rightarrow\infty
  \end{equation*}
  where $B(.,.)$ denotes the beta function.
\end{theo}

\vspace{0.5cm}

\Proof Let $k\in\mathbb{N}\backslash\{0\}$ and $\alpha\in[1,2)$ be
fixed. Since $\mu_1<\infty$, it follows that
$\varphi^{(1)}(0)=-\mu_1$ and $1-\varphi(s)\underset{s\downarrow
0}{\sim}\mu_1 \, s$. Hence,
$s=1-\varphi(\psi(s))\underset{s\downarrow 0}{\sim} \mu_1 \,
\psi(s)$ and we deduce that $\psi(s)\underset{s\downarrow 0}{\sim}
\frac{s}{\mu_1} \in \mathrm{RV}_1^0$. Obviously, we have
$\varphi^{(1)}(s)\underset{s\downarrow 0}{\sim} -\mu_1 \in
\mathrm{RV}_0^0$. Consequently, we get $\varphi^{(1)}\circ \psi
\in \mathrm{RV}_0^0$ and $\varphi^{(n)}\circ\psi \in
\mathrm{RV}_{\alpha-n}^0$ for $n>\alpha$, since $\varphi^{(n)} \in
\mathrm{RV}_{\alpha-n}^0$ by Lemma \ref{lemmaderiv}, $\psi\in
\mathrm{RV}_1^0$ and $\lim_{s\rightarrow 0}\psi(s)=0$.\\

Using Potter's theorem, we therefore obtain the following upper
bound for the integrand in $B_t(k_1, \ldots, k_r)$. Set
$\delta_r:=\frac{\zeta}{2k+r}$ with $\zeta=\epsilon$ if
$\epsilon\in(0,1)$ or $\zeta\in(0,1)$ otherwise. For this chosen
$\delta_r>0$, there exists $C_r=C_r(\delta_r)>1$ such that for all
$t>0$:
\begin{eqnarray*}
  g_t(w;k_1, \ldots, k_r) &\leq& C_r \, w^{2k-1}\left(\max\left\{w^{\delta_r},
  w^{-\delta_r}\right\}\right)^{2k-1} w^{r\alpha-2k}\left(\max\left\{w^{\delta_r},
  w^{-\delta_r}\right\}\right)^{r+1} q_r(w) \cr
  &=& C_r \, w^{r\alpha-1} \, \max\left\{w^{\zeta}, w^{-\zeta}\right\} \, q_r(w)=:h(w).
\end{eqnarray*}

Now, $\int_0^\infty h(w) \, dw<\infty$ if and only if $\int_0^1
w^{r\alpha-1-\zeta} \, q_r(w) \, dw<\infty$ and $\int_1^\infty
w^{r\alpha-1+\zeta} \, q_r(w) \, dw<\infty$.\\

Since $\zeta\in(0,1)$ and $\zeta\leq \epsilon$, we use
(\ref{intbeta}) together with Lemma \ref{lemmamomcond} to get:
\begin{equation*}
  \int_0^1 w^{r\alpha-1-\zeta} \, q_r(w) \, dw \leq
  \int_0^1 w^{r-1-\zeta} \, q_r(w) \, dw \leq \int_0^\infty w^{r-1-\zeta} \, q_r(w) \, dw
  =\Gamma(r-\zeta) \, \mathbb{E}\left\{\Lambda^{\zeta}\right\} < \infty
\end{equation*}
and since $-k(\alpha-1)-\epsilon \leq -r(\alpha-1)-\zeta <0$:
\begin{equation*}
  \int_1^\infty w^{r\alpha-1+\zeta} \, q_r(w) \, dw \leq \int_0^\infty w^{r\alpha-1+\zeta} \, q_r(w) \, dw
  =\Gamma(r\alpha+\zeta) \, \mathbb{E}\left\{\Lambda^{-r(\alpha-1)-\zeta}\right\} < \infty.
\end{equation*}

Hence, the function $h$ is integrable.\\

Finally, $\lim_{t\rightarrow\infty}g_t(w;k_1, \ldots, k_r)=
w^{2k-1}w^{r\alpha-2k}q_r(w)=w^{r\alpha-1}q_r(w)$. Thus, applying
Lebesgue's theorem on dominated convergence and using
(\ref{intbeta}), we deduce:
\begin{equation*}
  \lim_{t\rightarrow\infty}\int_0^\infty g_t(w;k_1, \ldots, k_r) \, dw
  =\int_0^\infty w^{r\alpha-1} \, q_r(w) \, dw
  =\Gamma(r\alpha) \, \mathbb{E}\left\{\Lambda^{r(1-\alpha)}\right\}.
\end{equation*}

Since $\ell\left(\frac{1}{s}\right)\in\mathrm{RV}_0^0$, the
uniform convergence theorem for slowly varying functions states
that $\ell\left(\frac{x}{s}\right)\underset{s\downarrow
0}{\sim}\ell\left(\frac{1}{s}\right)$ uniformly on each compact
$x$-set in $(0,\infty)$. Since $\lim_{s\rightarrow
0}\frac{s}{\psi(s)}=\mu_1\in(0,\infty)$, we consequently get
$\ell\left(\frac{1}{\psi(s)}\right)=\ell\left(\frac{s}{\psi(s)}\frac{1}{s}\right)
\underset{s\downarrow 0}{\sim}\ell\left(\frac{1}{s}\right)$. This
together with Lemma \ref{lemmaderiv} and
$\varphi^{(1)}(\psi(s))\underset{s\downarrow 0}{\sim} -\mu_1$
yields:
\begin{eqnarray*}
  f_t(k_1, \ldots, k_r) &\underset{t\uparrow\infty}{\sim}& \frac{\alpha^r}{\mu_1}
  \prod_{i=1}^{r}\Gamma(2k_i-\alpha) \, t^{r-1} \, \psi^{r\alpha-1}\left(\frac{1}{t}\right)
  \ell^r\left(\frac{1}{\psi\left(\frac{1}{t}\right)}\right) \cr
  &\underset{t\uparrow\infty}{\sim}& \frac{\alpha^r}{\mu_1^{r\alpha}}
  \prod_{i=1}^{r}\Gamma(2k_i-\alpha) \, t^{r(1-\alpha)} \, \ell^r(t).
\end{eqnarray*}

Therefore, we obtain:
\begin{equation*}
  B_t(k_1, \ldots, k_r) \underset{t\uparrow\infty}{\sim}
  \left(\frac{\alpha}{\mu_1^{\alpha}}\right)^r
  \Gamma(r\alpha) \, \mathbb{E}\left\{\Lambda^{r(1-\alpha)}\right\}
  \prod_{i=1}^{r}\Gamma(2k_i-\alpha) \, t^{r(1-\alpha)} \, \ell^r(t).
\end{equation*}

When $\alpha=1$, we have $\ell(t)=o(1)$ since $\mu_1<\infty$.
Hence, the first order asymptotic behavior of (\ref{momT}) is
solely determined by the term with $r=1$ and we obtain:
\begin{equation*}
  \mathbb{E}\left\{T_{N(t)}^k\right\} \sim
  \frac{\alpha}{\mu_1^\alpha}\frac{\Gamma(2k-\alpha)\Gamma(\alpha)}{(2k-1)!}
  \, \mathbb{E}\left\{\Lambda^{1-\alpha}\right\} \, t^{1-\alpha} \, \ell(t)
  \quad \text{as} \quad t\rightarrow\infty.
\end{equation*}

Finally, note that
$\mathbb{E}\left\{\Lambda^{1-\alpha}\right\}<\infty$ since
$-k(\alpha-1)-\epsilon<1-\alpha<0$ if $\alpha\neq 1$ and
$1-\alpha=0$ if $\alpha=1$. \ $\square$

\vspace{0.5cm}

We pass to the case $\alpha>2$.

\vspace{0.5cm}

\begin{theo}\label{theo>2}
  Assume that $X_1$ is of Pareto-type with index $\alpha>2$. Let $\{N(t); \, t\geq 0\}$ be
  a mixed Poisson process with mixing random variable $\Lambda$. Let $k\in\mathbb{N}\backslash\{0\}$
  be fixed. If $k=1$, assume further that $\mathbb{E}\left\{\Lambda^{-1-\epsilon}\right\}<\infty$ for some
  $\epsilon>0$. If $k\neq 1$, assume further that $\mathbb{E}\left\{\Lambda^{k-2+\epsilon}\right\}<\infty$
  and $\mathbb{E}\left\{\Lambda^{1-2k-\epsilon}\right\}<\infty$ for some $\epsilon>0$. Then for $k<\alpha-1$:
  \begin{equation}\label{k<alpha-1}
    \mathbb{E}\left\{T_{N(t)}^k\right\} \sim \left(\frac{\mu_2}{\mu_1^2}\right)^k
    \mathbb{E}\left\{\Lambda^{-k}\right\} \, t^{-k} \quad \text{as} \quad t\rightarrow\infty
  \end{equation}
  and for $k>\alpha-1$:
  \begin{equation}\label{k>alpha-1}
    \mathbb{E}\left\{T_{N(t)}^k\right\}\sim
    \frac{\alpha}{\mu_1^\alpha} \, B(2k-\alpha, \alpha) \, \mathbb{E}\left\{\Lambda^{1-\alpha}\right\}
    \, t^{1-\alpha} \, \ell(t) \quad \text{as} \quad t\rightarrow\infty.
  \end{equation}
  If $k=\alpha-1$, then:
  \begin{itemize}
    \item[$(i)$] (\ref{k<alpha-1}) holds if $\ell(x)=o(1)$ (and in particular if $\mu_{k+1}<\infty$);
    \item[$(ii)$] $\mathbb{E}\left\{T_{N(t)}^k\right\}\sim
    \left(\left(\frac{\mu_2}{\mu_1^2}\right)^k + C \, \frac{(k+1) \, B(k-1, k+1)}{\mu_1^{k+1}}\right)
    \mathbb{E}\left\{\Lambda^{-k}\right\} t^{-k}$ \text{as} $t\rightarrow\infty$
    holds if $\lim_{x\rightarrow\infty}\ell(x)=C$ for a positive constant $C$;
    \item[$(iii)$] (\ref{k>alpha-1}) holds otherwise.
  \end{itemize}
\end{theo}

\vspace{0.5cm}

\Proof Let $k\in\mathbb{N}\backslash\{0\}$ and $\alpha>2$ be
fixed. Since $\mu_1<\infty$, it follows that
$\varphi^{(1)}(0)=-\mu_1$ and $1-\varphi(s)\underset{s\downarrow
0}{\sim}\mu_1 \, s$. Hence,
$s=1-\varphi(\psi(s))\underset{s\downarrow 0}{\sim} \mu_1 \,
\psi(s)$ and we deduce that $\psi(s)\underset{s\downarrow 0}{\sim}
\frac{s}{\mu_1} \in \mathrm{RV}_1^0$, leading to
$\lim_{s\rightarrow 0}\psi(s)=0$. Since
$\varphi^{(1)}(s)\underset{s\downarrow 0}{\sim} -\mu_1 \in
\mathrm{RV}_0^0$, we easily get $\varphi^{(1)}\circ \psi \in
\mathrm{RV}_0^0$. Moreover, by Lemma \ref{lemmaderiv},
$\varphi^{(n)} \in \mathrm{RV}_{\alpha-n}^0$ for $n>\alpha$ and
$\varphi^{(n)} \in \mathrm{RV}_0^0$ for $n\leq\alpha$.
Consequently, we get $\varphi^{(n)}\circ\psi \in
\mathrm{RV}_{\alpha-n}^0$ for $n>\alpha$ and
$\varphi^{(n)}\circ\psi \in \mathrm{RV}_0^0$ for $n\leq\alpha$.\\

For simplicity, we first assume that $\alpha\notin\mathbb{N}$.
Using Potter's theorem, we obtain the following upper bound for
the integrand in $B_t(k_1, \ldots, k_r)$. Set
$\delta_r:=\frac{\zeta}{2k+r}$ with $\zeta=\epsilon$ if
$\epsilon\in(0,1)$ or $\zeta\in(0,1)$ otherwise. For this chosen
$\delta_r>0$, there exists $C_r=C_r(\delta_r)>1$ such that for all
$t>0$:
\begin{eqnarray*}
  g_t(w;k_1, \ldots, k_r) &\leq& C_r \, w^{2k-1}\left(\max\left\{w^{\delta_r},
  w^{-\delta_r}\right\}\right)^{2k-1} w^{r_1\alpha-2u_1}\left(\max\left\{w^{\delta_r},
  w^{-\delta_r}\right\}\right)^{r+1} q_r(w) \cr
  &=& C_r \, w^{2k-1+r_1\alpha-2u_1} \, \max\left\{w^{\zeta}, w^{-\zeta}\right\} \, q_r(w) =:h(w)
\end{eqnarray*}
where $r_1$ denotes the number of integers among $\{k_1, \ldots,
k_r\}$ that are greater than $\alpha/2$ and $u_1$ is the sum of
these.\\

Now, $\int_0^\infty h(w)\, dw<\infty$ if and only if $\int_0^1
w^{2k-1+r_1\alpha-2u_1-\zeta} \, q_r(w)\, dw<\infty$ and
$\int_1^\infty w^{2k-1+r_1\alpha-2u_1+\zeta} \, q_r(w)\,
dw<\infty$.\\

We have $2-2k\leq r_1\alpha-2u_1\leq 0$. Indeed, if $r_1=0$ then
obviously $r_1\alpha-2u_1=0$. Now, if $r_1\neq 0$ (i.e. $r_1\geq
1$) then $r_1\alpha-2u_1<0$ on the one hand, and
$r_1\alpha-2u_1\geq \alpha-2k$ on the other hand. Consequently,
$1\leq 2k-1+r_1\alpha-2u_1\leq 2k-1$.\\

Since $\zeta\in(0,1)$ and $\zeta\leq \epsilon$, we get
$0<r-2+\zeta\leq k-2+\epsilon$ if $r\geq 2$ and $1-2k-\epsilon\leq
-1+\zeta<0$ if $r=1$. Moreover, $1-2k-\epsilon \leq r-2k-\zeta<0$.
Therefore, using (\ref{intbeta}) together with Lemma
\ref{lemmamomcond} leads to:
\begin{equation*}
  \int_0^1 w^{2k-1+r_1\alpha-2u_1-\zeta} \, q_r(w) \, dw
  \leq \int_0^1 w^{1-\zeta} \, q_r(w) \, dw
  \leq \Gamma(2-\zeta)\, \mathbb{E}\left\{\Lambda^{r-2+\zeta}\right\}<\infty
\end{equation*}
and:
\begin{equation*}
  \int_1^\infty w^{2k-1+r_1\alpha-2u_1+\zeta} \, q_r(w) \, dw
  \leq \int_1^\infty w^{2k-1+\zeta} \, q_r(w) \, dw
  \leq \Gamma(2k+\zeta)\, \mathbb{E}\left\{\Lambda^{r-2k-\zeta}\right\}<\infty.
\end{equation*}

When $k=1$, obviously $r=1$ and we get $-1-\epsilon\leq -1-\zeta <
-1+\zeta<0$. The condition
$\mathbb{E}\left\{\Lambda^{-1-\epsilon}\right\}<\infty$ is thus
sufficient for
$\mathbb{E}\left\{\Lambda^{-1+\zeta}\right\}<\infty$ and
$\mathbb{E}\left\{\Lambda^{-1-\zeta}\right\}<\infty$.\\

Hence, the function $h$ is integrable.\\

Finally, $\lim_{t\rightarrow\infty}g_t(w;k_1, \ldots, k_r)=
w^{2k-1}w^{r_1\alpha-2u_1}q_r(w)=w^{2(k-u_1)+r_1\alpha-1}q_r(w)$.
Thus, applying Lebesgue's theorem on dominated convergence and
using (\ref{intbeta}), we deduce:
\begin{equation*}
  \lim_{t\rightarrow\infty}\int_0^\infty g_t(w;k_1, \ldots, k_r) \, dw
  =\int_0^\infty w^{2(k-u_1)+r_1\alpha-1} q_r(w) \, dw
  =\Gamma(2(k-u_1)+r_1\alpha) \, \mathbb{E}\left\{\Lambda^{r-2(k-u_1)-r_1\alpha}\right\}.
\end{equation*}

By virtue of the uniform convergence theorem for slowly varying
functions, we get $\ell\left(1/\psi(s)\right)\underset{s\downarrow
0}{\sim}\ell\left(\frac{1}{s}\right)$. This together with Lemma
\ref{lemmaderiv} and $\varphi^{(1)}(\psi(s))\underset{s\downarrow
0}{\sim} -\mu_1$ yields:
\begin{eqnarray}\label{alphanotinN}
  f_t(k_1, \ldots, k_r) &\underset{t\uparrow\infty}{\sim}& \frac{\alpha^{r_1}K_1
  K_2}{\mu_1} \, t^{r-1} \, \psi^{2(k-u_1)+r_1\alpha-1}\left(\frac{1}{t}\right)
  \ell^{r_1}\left(\frac{1}{\psi\left(\frac{1}{t}\right)}\right) \cr
  &\underset{t\uparrow\infty}{\sim}& \frac{\alpha^{r_1}K_1
  K_2}{\mu_1^{2(k-u_1)+r_1\alpha}} \, t^{r-2(k-u_1)-r_1\alpha} \, \ell^{r_1}(t)
\end{eqnarray}
where $K_1:=\prod_{j\,\in\,\{i: \,
2k_i>\alpha\}}\Gamma(2k_j-\alpha)$ and $K_2:=\prod_{j\,\in\,\{i:
\, 2k_i<\alpha\}}\mu_{2k_j}$, with $\text{card}\{i: \,
2k_i>\alpha\}=r_1$ and $\text{card}\{i: \, 2k_i<\alpha\}=r-r_1$.
By convention $\prod_{j\,\in\,\emptyset}=1$ and
$\text{card}\,\emptyset=0$.\\

Therefore, the asymptotic behavior of $B_t(k_1, \ldots, k_r)$ is
given by:
\begin{equation*}
  B_t(k_1, \ldots, k_r) \underset{t\uparrow\infty}{\sim}
  \frac{\alpha^{r_1}K_1 K_2}{\mu_1^{2(k-u_1)+r_1\alpha}} \,
  \Gamma(2(k-u_1)+r_1\alpha) \,
  \mathbb{E}\left\{\Lambda^{r-2(k-u_1)-r_1\alpha}\right\} \, t^{r-2(k-u_1)-r_1\alpha} \,
  \ell^{r_1}(t).
\end{equation*}

It remains to determine the dominating asymptotic term among all
possible $B_t(k_1, \ldots, k_r)$. For $r_1>0$, the largest
exponent is achieved with $r_1=1$, $u_1=k$ and thus $r=1$, so that
the asymptotic order is $t^{1-\alpha}\ell(t)$. Note that $r_1>0$
is possible for $2k>\alpha$ only. For $r_1=0$, obviously $r=k$
(which implies $k_1=\cdots=k_r=1$) dominates, leading to the
asymptotic order $t^{-k}$. Hence, the asymptotically dominating
power among all $B_t(k_1, \ldots, k_r)$ is given by
$\max\{1-\alpha, -k\}$. From this, we see that when $k<\alpha-1$,
$r=k$ dominates and we obtain from (\ref{momT}):
\begin{equation*}
  \mathbb{E}\left\{T_{N(t)}^k\right\} \sim
  \frac{k! \, \mu_2^k \, \Gamma(2k) \, \mathbb{E}\left\{\Lambda^{-k}\right\}}
  {(2k-1)! \, k! \, \mu_1^{2k}} \, t^{-k}
  = \left(\frac{\mu_2}{\mu_1^2}\right)^k \mathbb{E}\left\{\Lambda^{-k}\right\} \, t^{-k}
  \quad \text{as} \quad t\rightarrow\infty.
\end{equation*}

Alternatively, when $k>\alpha-1$, the term with $r=1$ dominates
and we find:
\begin{equation*}
  \mathbb{E}\left\{T_{N(t)}^k\right\}\sim
  \frac{k! \, \alpha \, \Gamma(2k-\alpha) \, \Gamma(\alpha) \,
  \mathbb{E}\left\{\Lambda^{1-\alpha}\right\}}{k! \, (2k-1)! \, \mu_1^\alpha} \, t^{1-\alpha} \, \ell(t)
  = \frac{\alpha}{\mu_1^\alpha} \, B(2k-\alpha, \alpha) \, \mathbb{E}\left\{\Lambda^{1-\alpha}\right\}
  \, t^{1-\alpha} \, \ell(t) \quad \text{as} \quad t\rightarrow\infty
\end{equation*}
which is the same expression as the one obtained in Theorem
\ref{theo12} for $\alpha\in[1,2)$ and $\mu_1<\infty$.\\

The above conclusions also hold for $\alpha\in\mathbb{N}$ as long
as $k\neq\alpha-1$. Nevertheless, just note that instead of
(\ref{alphanotinN}) we have the following by virtue of Lemma
\ref{lemmaderiv}:
\begin{equation*}
  f_t(k_1, \ldots, k_r) \underset{t\uparrow\infty}{\sim} \frac{\alpha^{r_1+r_2}K_1
  K_3}{\mu_1^{2(k-u_1)+r_1\alpha}} \, t^{r-2(k-u_1)-r_1\alpha} \, \ell^{r_1}(t) \, \tilde{\ell}^{r_2}(t)
\end{equation*}
where $\tilde{\ell}(x)=\int_0^x \frac{\ell(u)}{u} \, du \in
\mathrm{RV}_0^\infty$, $K_3:=\prod_{j\,\in\,\{i: \,
2k_i\leq\alpha, \, \mu_{2k_i}<\infty\}}\mu_{2k_j}$ and
$r_2:=\text{card}\{i: \, 2k_i=\alpha, \, \mu_{2k_i}=\infty\}$,
with $\text{card}\{i: \, 2k_i\leq\alpha, \,
\mu_{2k_i}<\infty\}=r-r_1-r_2$.\\

When $k=\alpha-1$, the slowly varying function $\ell$ determines
which of the two terms $t^{1-\alpha}\ell(t)$ (corresponding to
$r=1$) and $t^{-k}$ (corresponding to $r=k$) dominates the
asymptotic behavior. If $\ell(x)=o(1)$, which is in particular
fulfilled if $\mu_{k+1}<\infty$, the second term dominates. If
$\lim_{x\rightarrow\infty}\ell(x)=C$ for a positive constant $C$,
then both terms matter. Otherwise, the first term dominates.\\

To end the proof, it remains to check that
$\mathbb{E}\left\{\Lambda^{-k}\right\}<\infty$ for $k\leq\alpha-1$
and that $\mathbb{E}\left\{\Lambda^{1-\alpha}\right\}<\infty$ for
$k>\alpha-1$. When $k>\alpha-1$, we have
$1-2k-\epsilon<-k<1-\alpha<0$. Otherwise, we have
$1-2k-\epsilon<-k<0$. Thus, we conclude by using Lemma
\ref{lemmamomcond}. \ $\square$

\vspace{0.5cm}

Finally, we consider the case $\alpha=2$.

\vspace{0.5cm}

\begin{coro}\label{coro2}
  Assume that $X_1$ is of Pareto-type with index $\alpha=2$. Let $\{N(t); \, t\geq 0\}$ be a
  mixed Poisson process with mixing random variable $\Lambda$.
  \begin{itemize}
    \item[$(i)$] If $\mathbb{E}\left\{\Lambda^{-1-\epsilon}\right\}<\infty$ for some $\epsilon>0$, then:
    \begin{equation*}
      \mathbb{E}\left\{T_{N(t)}\right\} \sim
      \begin{cases}
        \frac{\mu_2 \, \mathbb{E}\left\{\Lambda^{-1}\right\}}{\mu_1^2} \frac{1}{t}
        & \text{if $\mu_2<\infty$}\\
        \frac{2 \, \mathbb{E}\left\{\Lambda^{-1}\right\}}{\mu_1^2} \frac{\tilde{\ell}(t)}{t}
        & \text{if $\mu_2=\infty$}
      \end{cases}
      \quad \text{as} \quad t\rightarrow\infty
    \end{equation*}
    where $\tilde{\ell}(x)=\int_0^x \frac{\ell(u)}{u} \, du \in \mathrm{RV}_0^\infty$.
    \item[$(ii)$] For any fixed $k\geq 2$, if $\mathbb{E}\left\{\Lambda^{k-2+\epsilon}\right\}<\infty$
    and $\mathbb{E}\left\{\Lambda^{1-2k-\epsilon}\right\}<\infty$ for some $\epsilon>0$, then:
    \begin{equation*}
      \mathbb{E}\left\{T_{N(t)}^k\right\} \sim \frac{\mathbb{E}\left\{\Lambda^{-1}\right\}}
      {\mu_1^2 \, (k-1) \, (2k-1)} \frac{\ell(t)}{t} \quad \text{as} \quad t\rightarrow\infty.
    \end{equation*}
  \end{itemize}
\end{coro}

\vspace{0.5cm}

\Proof One can easily verify that Theorem \ref{theo>2} remains
true for $\alpha=2$, except when $k=1$ if $\mu_2=\infty$. In the
latter case, obviously $r=1$ and $B_t(k_1, \ldots, k_r)$ becomes:
\begin{eqnarray*}
  B_t(1) &=& \int_0^\infty s \, \varphi^{(2)}(s) \, Q_t^{(1)}(\varphi(s)) \, ds \cr
  &=& \underbrace{\frac{-\psi\left(\frac{1}{t}\right)}
  {\varphi^{(1)}\left(\psi\left(\frac{1}{t}\right)\right)}
  \, \varphi^{(2)}\left(\psi\left(\frac{1}{t}\right)\right)}_{=:f_t(1)}
  \int_0^\infty \underbrace{\frac{\psi\left(\frac{w}{t}\right)}{\psi\left(\frac{1}{t}\right)}
  \, \frac{\varphi^{(2)}\left(\psi\left(\frac{w}{t}\right)\right)}
  {\varphi^{(2)}\left(\psi\left(\frac{1}{t}\right)\right)}
  \, \frac{\varphi^{(1)}\left(\psi\left(\frac{1}{t}\right)\right)}
  {\varphi^{(1)}\left(\psi(\frac{w}{t})\right)} \, q_1(w) \, \mathbbm{1}_{[0,t)}(w)}_{=:g_t(w;1)} \, dw.
\end{eqnarray*}

Since $\mu_1<\infty$, it follows that $\varphi^{(1)}(0)=-\mu_1$
and $1-\varphi(s)\underset{s\downarrow 0}{\sim}\mu_1 \, s$. Hence,
$s=1-\varphi(\psi(s))\underset{s\downarrow 0}{\sim} \mu_1 \,
\psi(s)$ and we deduce that $\psi(s)\underset{s\downarrow 0}{\sim}
\frac{s}{\mu_1} \in \mathrm{RV}_1^0$, leading to
$\lim_{s\rightarrow 0}\psi(s)=0$. Since
$\varphi^{(1)}(s)\underset{s\downarrow 0}{\sim} -\mu_1 \in
\mathrm{RV}_0^0$, we easily get $\varphi^{(1)}\circ \psi \in
\mathrm{RV}_0^0$. Moreover, $\varphi^{(2)}(s)\underset{s\downarrow
0}{\sim} 2 \, \tilde{\ell}\left(\frac{1}{s}\right) \in
\mathrm{RV}_0^0$ by Lemma \ref{lemmaderiv} and as a consequence
$\varphi^{(2)}\circ \psi \in \mathrm{RV}_0^0$, where
$\tilde{\ell}(x)=\int_0^x
\frac{\ell(u)}{u} \, du \in \mathrm{RV}_0^\infty$.\\

Using Potter's theorem, we therefore obtain the following upper
bound for the integrand in $B_1(1)$. Set $\delta:=\frac{\zeta}{3}$
with $\zeta=\epsilon$ if $\epsilon\in(0,1)$ or $\zeta\in(0,1)$
otherwise. For this chosen $\delta>0$, there exists
$C=C(\delta)>1$ such that for all $t>0$:
\begin{eqnarray*}
  g_t(w;1) &\leq& C \, w \left(\max\left\{w^{\delta}, w^{-\delta}\right\}\right)^3 q_1(w) \cr
  &=& C \, w \, \max\left\{w^{\zeta}, w^{-\zeta}\right\} \, q_1(w) =:h(w).
\end{eqnarray*}

Now, $\int_0^\infty h(w)\, dw<\infty$ if and only if $\int_0^1
w^{1-\zeta} \, q_1(w)\, dw<\infty$ and $\int_1^\infty w^{1+\zeta}
\, q_1(w)\, dw<\infty$.\\

Since $\zeta\in(0,1)$ and $\zeta\leq \epsilon$, we get
$-1-\epsilon \leq -1-\zeta < -1+\zeta <0$. Therefore, using
(\ref{intbeta}) together with Lemma \ref{lemmamomcond} leads to:
\begin{equation*}
  \int_0^1 w^{1-\zeta} \, q_1(w) \, dw \leq \int_0^\infty w^{1-\zeta} \, q_1(w)\, dw
  =\Gamma(2-\zeta)\, \mathbb{E}\left\{\Lambda^{-1+\zeta}\right\}<\infty
\end{equation*}
and:
\begin{equation*}
  \int_1^\infty w^{1+\zeta} \, q_1(w) \, dw \leq \int_0^\infty w^{1+\zeta} \, q_1(w)\, dw
  =\Gamma(2+\zeta)\, \mathbb{E}\left\{\Lambda^{-1-\zeta}\right\}<\infty.
\end{equation*}

Hence, the function $h$ is integrable.\\

Finally, $\lim_{t\rightarrow\infty}g_t(w;1)= w \, q_1(w)$. Thus,
applying Lebesgue's theorem on dominated convergence and using
(\ref{intbeta}), we deduce:
\begin{equation*}
  \lim_{t\rightarrow\infty}\int_0^\infty g_t(w;1) \, dw
  =\int_0^\infty w \, q_1(w) \, dw=\mathbb{E}\left\{\Lambda^{-1}\right\}.
\end{equation*}

By virtue of the uniform convergence theorem for slowly varying
functions, we get
$\tilde{\ell}\left(1/\psi(s)\right)\underset{s\downarrow
0}{\sim}\tilde{\ell}\left(\frac{1}{s}\right)$. This together with
$\varphi^{(1)}(\psi(s))\underset{s\downarrow 0}{\sim} -\mu_1$ and
$\varphi^{(2)}(\psi(s))\underset{s\downarrow 0}{\sim} 2 \,
\tilde{\ell}\left(1/\psi(s)\right)$ yields:
\begin{equation*}
  f_t(1) \underset{t\uparrow\infty}{\sim} \frac{2}{\mu_1^2}\frac{\tilde{\ell}(t)}{t}.
\end{equation*}

Consequently, we obtain:
\begin{equation*}
  \mathbb{E}\left\{T_{N(t)}\right\} = B_t(1) \sim
  \frac{2 \, \mathbb{E}\left\{\Lambda^{-1}\right\}}{\mu_1^2}\frac{\tilde{\ell}(t)}{t}
  \quad \text{as} \quad t\rightarrow\infty.
\end{equation*}

Finally, note that $\mathbb{E}\left\{\Lambda^{-1}\right\}<\infty$
since $-1-\epsilon<-1$. \ $\square$

\vspace{0.5cm}

We end by a remark.

\paragraph{Remark} \textit{As in Albrecher and Teugels \cite{at04}, the
integral representation approach that we use in this paper does
not permit to get a general asymptotic result for
$\mathbb{E}\left\{T_{N(t)}^k\right\}$ when $F$ in the additive
domain of attraction of a normal law, i.e. when $F$ has a slowly
varying truncated variance. Note that if the distribution function
$F$ of $X_1$ is as in Corollary \ref{coro2}, then $F$ is in the
additive domain of attraction of a normal law.}

\section{Conclusion}\label{ccl}

In this paper, we have derived the limiting behavior of arbitrary
moments of $T_{N(t)}$ assuming the distribution function $F$ of
$X_1$ to be of Pareto-type with index $\alpha>0$ and the counting
process $\{N(t); \ t\geq 0\}$ to be mixed Poisson. Different
results have shown up depending on the range of $\alpha$.\\

In the special setting where the distribution function $H$ of the
mixing random variable $\Lambda$ is degenerate at the point 1, our
results are similar to those derived by Albrecher and Teugels
\cite{at04} who assume the counting process to be non-random. This
is basically explained by the fact that the important quantities
$q_r(w)$ are the same in these two cases.\\

The \textit{coefficient of variation} of a positive random
variable $X$ is defined and denoted by:
\begin{equation*}
  CoVar(X):=\frac{\sqrt{\mathbb{V}X}}{\mathbb{E}X}
\end{equation*}
where $\mathbb{V}X$ denotes the variance of $X$. This risk measure
is frequently used in practice and is very popular among
actuaries.\\

From a random sample $X_1, \ldots, X_{N(t)}$ from $X$ of random
size $N(t)$ from an integer valued distribution, the coefficient
of variation $CoVar(X)$ is naturally estimated by the sample
coefficient of variation of $X$ defined and denoted by:
\begin{equation*}
  \widehat{CoVar(X)}:=\frac{S}{\overline{X}}
\end{equation*}
where $\overline{X}:=\frac{X_1+\cdots+X_{N(t)}}{N(t)}$ is the
sample mean and
$S^2:=\frac{1}{N(t)}\sum_{i=1}^{N(t)}\left(X_i-\overline{X}\right)^2$
is the sample variance.\\

The properties of the sample coefficient of variation
$\widehat{CoVar(X)}$ are usually studied in assuming the
finiteness of sufficiently many moments of $X$. The existence of
moments of $X$ is not always guaranteed in practical applications.
Hence, it is useful to investigate asymptotic properties of
$\widehat{CoVar(X)}$ also in these cases and it turns out that
this can be done by using results on $T_{N(t)}$.\\

Indeed, it appears that the quantity $T_{N(t)}$ is a basic
ingredient in the study of the sample coefficient of variation
since the following holds:
\begin{equation}\label{relcov}
  \widehat{CoVar(X)}=\sqrt{N(t) \, T_{N(t)}-1}.
\end{equation}

In a forthcoming paper, we will touch on the question of
\textit{convergence in distribution} for the appropriately
normalized quantity $T_{N(t)}$ when $X$ is of Pareto-type with
positive index $\alpha$ and the counting process $\{N(t); \, t\geq
0\}$ is mixed Poisson. Thanks to relation (\ref{relcov}),
asymptotic properties of the sample coefficient of variation will
be derived, even when the first moment and/or the second moment of
$X$ do not exist.\\

Incidentally, it will also be seen how the methodology can be
adapted to derive asymptotic properties of another risk measure,
the \textit{sample dispersion}. Recall that the value of the
\textit{dispersion} allows to compare the volatility with respect
to the Poisson case.

\paragraph{Acknowledgments} The author is supported by the grant
BDB-B/04/03 of the Katholieke Universiteit Leuven. The author
would like to thank Jef Teugels and Krishanu Maulik for fruitful
discussions. Part of the work has been done at the research
institute EURANDOM (Eindhoven, The Netherlands).

\end{document}